\documentclass[12pt,a4paper]{article}
\usepackage{a4wide}
\usepackage{amsmath,amsthm,amsfonts,amssymb,bbm}
\usepackage{graphicx,psfrag}
\usepackage{cite}

\numberwithin{equation}{section}

\newtheorem{prop}{Proposition}
\newtheorem{thm}[prop]{Theorem}
\newtheorem{lem}[prop]{Lemma}

\newenvironment{proofOF}[2]{\removelastskip\vspace{6pt}\noindent {\it Proof of #1.}~\rm#2}{\qed \par\vspace{6pt}}

\title{Maximum of Dyson Brownian motion and non-colliding systems 
with a boundary}
\author{A. Borodin, P.L. Ferrari, M. Pr\"ahofer, T. Sasamoto and J. Warren}
\author{Alexei Borodin\thanks{California Institute of Technology, 
e-mail: borodin@caltech.edu},
Patrik L. Ferrari\thanks{Bonn University, 
e-mail: ferrari@uni-bonn.de},
Michael Pr\"ahofer \thanks{TU M\"unchen,
e-mail: praehofer@ma.tum.de},
Tomohiro Sasamoto\thanks{TU M\"unchen, e-mail: sasamoto@ma.tum.de, 
Chiba University, sasamoto@math.s.chiba-u.ac.jp},\\
Jon Warren\thanks{University of Warwick,
e-mail: j.warren@warwick.ac.uk}}

\date{May 25, 2009}

\begin{document}
\maketitle \sloppy

\begin{abstract}
We prove an equality-in-law relating the maximum of GUE Dyson's 
Brownian motion and the non-colliding systems with a wall. 
This generalizes the well known relation between the maximum of 
a Brownian motion and a reflected Brownian motion. 
\end{abstract}

\section{Introduction and Results}\label{intro}
Dyson's Brownian motion model of GUE (Gaussian unitary ensemble)
is a stochastic process of positions of $m$ particles, 
$X(t) = (X_1(t),\ldots,X_m(t))$ described by the stochastic 
differential equation, 
\begin{equation}
 dX_i = dB_i + 
        \sum_{\substack{1\leq j\leq m \\j\neq i}}
        \frac{dt}{X_i-X_j},\quad 1\leq i\leq m,
\end{equation}
where $B_i,1\leq i \leq m$ are independent one dimensional Brownian 
motions\cite{Dyson1962}. The process satisfies  $X_1(t)< X_2(t) <
\cdots <X_m(t)$ for all $t>0$. 
We remark that the process $X$ can be started from the origin, i.e.,
one can take $X_i(0)=0,1\leq i \leq m$. See \cite{OY2002}.

One can introduce similar non-colliding system of $m$ particles
with a  wall at the origin \cite{grabiner, KT2004,TW2007}.
The dynamics of the positions of the $m$ particles 
$X^{(C)}=(X_1^{(C)},\ldots,X_m^{(C)})$ satisfying $0<X_1(t)< X_2(t) <
\cdots <X_m(t)$ for all $t>0$ 
are described by the stochastic differential equation,  
\begin{equation}
 dX_i^{(C)} 
 = 
 dB_i + \frac{dt}{X_i^{(C)}} +\sum_{\substack{1\leq j \leq m \\j\neq i}}
 \left(\frac{1}{X_i^{(C)}-X_j^{(C)}}+\frac{1}{X_i^{(C)}+X_j^{(C)}}\right)dt, 
 ~1\leq i\leq m.
\end{equation}
This process is referred to as Dyson's Brownian motion of type $C$. It
can be interpreted as a system of $m$ Brownian  particles conditioned
to never collide with each other or the wall.
 
One can also consider the case where the  wall above 
is replaced by a reflecting wall\cite{KT2004}. 
The dynamics of the positions of the $m$ particles 
$X^{(D)}=(X_1^{(D)},\ldots,X_m^{(D)})$ satisfying $0\leq X_1(t)< X_2(t) <
\cdots <X_m(t)$ for all $t>0$, 
is described by the stochastic differential equation,  
\begin{equation}
 dX_i^{(D)} 
 = 
 dB_i + \frac12{\mathbf 1}_{(i=1)} dL(t) 
 +
 \sum_{\substack{1\leq j \leq m\\j\neq i}}
 \left(\frac{1}{X_i^{(D)}-X_j^{(D)}}
 +
 \frac{1}{X_i^{(D)}+X_j^{(D)}}\right)dt, ~1\leq i\leq m,
\end{equation}
where $L(t)$ denotes the local time of $X_1^{(D)}$ at the origin.
This process will be  referred to as Dyson's Brownian motion of type $D$.
Some authors consider a process defined by the s.d.e.s (1.3) without the 
local time term. In this case the first component of the process is not 
constrained to remain non-negative, and the process takes values in the 
Weyl chamber of type $D$, $\{ |x_1|< x_2< x_3 \ldots < x_m\}$. The process we 
consider with a reflecting wall is obtained from this by replacing the 
first component with its absolute value, with the local time term 
appearing as a consequence of Tanaka's formula.

It is known  the processes
$X^{(C,D)}$ can be obtained using the Doob $h$-transform, see
\cite{grabiner}.  Let $(P_t^{0,(C,D)};t\geq 0)$ be the transition 
semigroup for $m$ independent Brownian motions killed on exiting 
$\{0<x_1<x_2 \ldots <x_m\}$, resp. the transition 
semigroup for $m$ independent Brownian motions reflected at the origin 
killed on exiting $\{0\leq x_1<x_2 \ldots <x_m\}$. 
From the Karlin-McGregor formula, the corresponding densities can be
written as 
\begin{equation}
 \det\{\phi_t(x_i-x_j')-\phi_t(x_i+x_j')\}_{1\leq i,j \leq m},
\end{equation}
resp., 
\begin{equation}
 \det\{\phi_t(x_i-x_j')+\phi_t(x_i+x_j')\}_{1\leq i,j \leq m},
\end{equation}
where $\phi_t(z) = \frac{1}{\sqrt{2\pi t}}e^{-z^2/(2t)}$.
Let 
\begin{equation}
\begin{aligned}
         h^{(C)}(x)& = 
         \prod_{i=1}^m x_i\prod_{1\leq i<j \leq m} (x_j^2-x_i^2), \\
         h^{(D)}(x)& =
         \prod_{1\leq i<j \leq m} (x_j^2-x_i^2). 
\end{aligned}
\label{h}
\end{equation}
For notational simplicity we suppress the index $C,D$ for the semigroups
and in $h$ in the following. 
Then one can show that $h(x)$ is invariant for the $P_t^0$ semigroup 
and we may define a Markov semigroup by  
\begin{equation}\label{Pt}
 P_t(x,dx') = h(x')P_t^0(x,dx')/h(x).
\end{equation}
This is the semigroup of the Dyson non-colliding system of Brownian 
motions of type $C$ and $D$. 
Similarly to the $X$ process, the processes $X^{(C)}$ and $X^{(D)}$ 
can also be started from the origin (see \cite{BBO2005} or use Lemma 4 
in \cite{KT2004} and apply the same arguments as in \cite{OY2002}). 

In GUE Dyson's Brownian motion of $n$ particles, let us take the 
initial conditions to be $X_i(0)=0,1\leq i \leq n$. The quantity 
we are interested in is the maximum of the position of the 
top particle for a finite duration of time, 
$\max_{0\leq s \leq t} X_n(s)$. In the sequel we write $\sup$
instead of $\max$ to conform with common usage in the literature. 
Let $m$ be the integer such that $n=2m$ when $n$ is even
and $n=2m-1$ when $n$ is odd. Consider the non-colliding 
systems of $X^{(C)},X^{(D)}$ of $m$ particles starting from the 
origin, $X_i^{(C,D)}(0)=0,1\leq i\leq m$. 

Our main result of this note is 

\begin{thm}
Let $X$ and $X^{(C)},X^{(D)}$ start from the origin. Then 
for each fixed $t\geq 0$, one has 
\begin{equation}
\sup_{0\leq s \leq t} X_n(s) \overset{d}{=}
\begin{cases}
 X_m^{(C)}(t), &\text{for } n=2m,\\
 X_m^{(D)}(t), &\text{for } n=2m-1. 
\end{cases}
\label{thm}
\end{equation}
\end{thm}

To prove the theorem we introduce two more processes $Z_j$ and $Y_j$.
In the $Z$ process, $Z_1 \leq Z_2 \leq \ldots \leq Z_n$, $Z_1$ is a 
Brownian motion and $Z_{j+1}$ is reflected by $Z_j$, $1\leq j \leq n-1$.   
Here the reflection means the Skorokhod construction to push $Z_{j+1}$
up from $Z_j$. More precisely,
\begin{equation}
\begin{aligned}
 Z_1(t) &= B_1(t), \\
 Z_j(t) &= \sup_{0\leq s \leq t}(Z_{j-1}(s)+B_j(t)-B_j(s)),~~
 2 \leq j \leq n,
\end{aligned}
\label{ZinB}
\end{equation}
where $B_i,1\leq i \leq n$ are independent Brownian motions, each
starting from $0$.
The process is the same as the process 
$(X_1^1(t),X_2^2(t),\ldots,X_n^n(t);t \geq 0)$ studied in section 4 
of \cite{Warren2007}. The representation (\ref{ZinB}) was given 
earlier in \cite{Baryshnikov2001}. 
In the $Y$ process, $0 \leq Y_1 \leq Y_2 \leq \ldots \leq Y_n$, 
the interactions among $Y_i$'s are the same as in the $Z$ process, 
i.e., $Y_{j+1}$ is reflected by $Y_j$, $1\leq j \leq n-1$, 
but $Y_1$ is now a Brownian motion reflected at the origin 
(again by Skorokhod construction). 
Similarly to (\ref{ZinB}),
\begin{equation}
\begin{aligned}
 Y_1(t) &= B_1(t)-\inf_{0\leq s \leq t}B_1(s)
         = \sup_{0\leq s\leq t}(B_1(t)-B_1(s)), \\
 Y_j(t) &= \sup_{0\leq s \leq t}(Y_{j-1}(s)+B_j(t)-B_j(s)),~~
 2 \leq j \leq n.
\end{aligned}
 \label{YinB}
\end{equation}

From the results in \cite{OY2002,BJ2002,Warren2007}, we know 
\begin{equation}
 (X_n(t); t\geq 0) \overset{d}{=} (Z_n(t); t\geq 0)
\label{XZ}
\end{equation}
and hence
\begin{equation}
 \sup_{0\leq s \leq t} X_n(s) \overset{d}{=} 
 \sup_{0\leq s \leq t} Z_n(s).
\end{equation}

In this note we show 
\begin{prop}\label{WXDC}
The following equalities in law hold between processes:
\begin{equation}
\begin{aligned}
 (Y_{2m}(t); t \geq 0)  &\overset{d}{=} (X_m^{(C)}(t); t \geq 0), \\
 (Y_{2m-1}(t);t \geq 0) &\overset{d}{=} (X_m^{(D)}(t); t \geq 0), 
\end{aligned}
\label{YX}
\end{equation} 
$m\in\mathbb{N}$. 
\end{prop}
\noindent
The proof of this proposition is given in Section \ref{XCD_Wr}.
The idea behind it  is that the processes $(Y_i)_{i\ge 1}$ and 
$(X_j^{(C,D)})_{j\ge 1}$ could be realized on a common probability 
space consisting of Brownian motions satisfying certain interlacing 
conditions with a boundary \cite{Warren2007,WW2008}. Such a system is expected to 
appear as a  scaling limit of the discrete processes considered 
in \cite{BK2009,WW2008}. In this enlarged process, the processes
$Y_n(t)$ and $X_m^{(C,D)}(t)$ just represent two different ways of 
looking at the evolution of a specific particle and so  
the statement of Proposition 2 follows immediately. 
Justification of such an approach is however quite involved, and we 
prefer to give a simple independent proof. 
See also \cite{BJ2002} for another representation of $X_m^{(C,D)}$ in 
terms of independent Brownian motions.

Then to prove (\ref{thm}) it is enough to show
\begin{prop}\label{maxWWr}
For each fixed $t$ we have 
\begin{equation}
 \sup_{0\leq s \leq t} Z_n(s) \overset{d}{=} Y_n(t).
\end{equation} 
\end{prop}
\noindent
This is shown in Section \ref{WWr}.
For $n=1$ case, this is well known from the Skorokhod construction 
of reflected Brownian motion \cite{RY1999}. 
The $n>1$ case can also be understood 
graphically by reversing time direction and the order of particles. 
This relation could also be established as a limiting case of the last 
passage percolation. In fact the identities in our theorem was first 
anticipated from the consideration of a diffusion scaling limit of 
the totally asymmetric exclusion process with 2 speeds 
\cite{BFS2009p}.

\medskip\noindent
{\it Acknowledgments.}

AB was partially supported by NSF grant DMS-0707163. 
TS thanks S. Grosskinsky and O. Zaboronski for inviting him to a 
workshop at University of Warwick, and N. O'Connell and H. Spohn 
for useful discussions and suggestions. 
His work was partially supported by the Grant-in-Aid for Young
Scientists (B), the Ministry of Education, Culture, Sports, Science
and Technology, Japan.

\section{Proof of proposition \ref{WXDC}}
\label{XCD_Wr}
In this section we prove the relation between $X^{(C,D)}$ and 
$Y$, (\ref{YX}). The following Lemma is a 
generalization of the Rogers-Pitman criterion 
\cite{RP1981} for a function of a Markov process to be Markovian. 

\begin{lem}\label{RP}
Suppose that $\{X(t):t\geq 0\}$ is a Markov process with state
space $E$, evolving according to a transition semigroup 
$(P_t;t\geq 0)$ and with initial distribution $\mu$. Suppose 
that $\{Y(t):t\geq 0\}$ is a Markov process with state space
$F$, evolving according to a transition semigroup $(Q_t;t\geq 0)$
and with initial distribution $\nu$. Suppose further that $L$
is a Markov transition kernel from $E$ to $F$, such that 
$\mu L=\nu$ and the intertwining $P_t L=L Q_t$ holds. Now let
$f: E\to G$ and $g: F\to G$ be maps into a third state space $G$,
and suppose that 
\begin{center}
 $L(x,\cdot)$ is carried by $\{y\in F: g(y)=f(x)\}$ for each $x\in E$.
\end{center}
Then we have 
\begin{center}
 $\{ f(X(t)):t\geq 0\} \overset{d}{=} \{g(Y(t)):t\geq 0\}$,
\end{center}
in the sense of finite dimensional distributions. 
\end{lem}

\begin{proofOF}{Lemma \ref{RP}}
For any bounded function $\alpha$ on $G$ let $\Gamma_1 \alpha$ be 
the function $\alpha \circ f$ defined on $E$ and let 
$\Gamma_2 \alpha$ be the function $\alpha \circ g$ defined on $F$.
Then it follows from the condition that $L(x,\cdot)$ is carried by
$\{y\in F: g(y)=f(x)\}$ that whenever $h$ is a bounded function defined
on $F$ then 
\begin{equation}
 L(\Gamma_2 \alpha \times h) = \Gamma_1 \alpha \times Lh,
\end{equation}
which is shorthand for 
$\int L(x,dy)\Gamma_2 \alpha(y) h(y) = \Gamma_1 \alpha \times Lh$.
For any bounded test functions $\alpha_0,\alpha_1,\cdots,\alpha_n$
defined on $G$, and times $0<t_1<\cdots <t_n,$ we have, using the 
previous equation and the intertwining relation repeatedly, 
\begin{align}
& \quad
\mathbb{E}[\alpha_0(g(Y(0))) \alpha_1(g(Y(t_1))) \ldots \alpha_n(g(Y(t_n)))]
\notag\\
&=
\nu(\Gamma_2\alpha_0 \times Q_{t_1}(\Gamma_2 \alpha_1 \times Q_{t_2-t_1}
   (\cdots (\Gamma_2 \alpha_{n-1}\times Q_{t_n-t_{n-1}}\Gamma_2\alpha_n)
    \cdots )))
\notag\\
&=
\mu L(\Gamma_2\alpha_0 \times Q_{t_1}(\Gamma_2 \alpha_1 \times Q_{t_2-t_1}
   (\cdots (\Gamma_2 \alpha_{n-1}\times Q_{t_n-t_{n-1}}\Gamma_2\alpha_n)
    \cdots )))
\notag\\
&=
\mu(\Gamma_1\alpha_0 \times P_{t_1}(\Gamma_1 \alpha_1 \times P_{t_2-t_1}
   (\cdots (\Gamma_1 \alpha_{n-1}\times P_{t_n-t_{n-1}}\Gamma_1\alpha_n)
    \cdots )))
\notag\\
&=
 \mathbb{E}[\alpha_0(f(X(0))) \alpha_1(f(X(t_1))) \ldots \alpha_n(f(X(t_n)))] 
\end{align}
which proves the equality in law.  
\end{proofOF}

We let $(Y(t):t\geq 0)$ be the process $Y$ of $n$ reflected Brownian 
motions with a wall introduced in the previous section. 
It is clear from the construction (\ref{YinB}) that the process $Y$
is a time homogeneous Markov process. We denote its transition
semigroup  by $\bigl(Q_t; t \geq 0)$. It turns out that there is an explicit
formula for the corresponding densities. Recall   
$\phi_t(z) = \frac{1}{\sqrt{2\pi t}}e^{-z^2/(2t)}$. Let us define 
$\phi_t^{(k)}(y) = \frac{d^k}{d y^k} \phi_t(y)$ for $k\geq 0$
and $\phi_t^{(-k)}(y) = (-1)^k\int_y^{\infty} 
\frac{(z-y)^{k-1}}{(k-1)!} \phi_t(z) dz$ for $k\geq 1$.

\begin{prop}\label{GY}
The transition densities $q_t(y,y')$ from $y=(y_1,\ldots,y_n)$
at $t=0$ to $y'=(y_1',\ldots,y_n')$ at $t$ of the $Y$ process can 
be written as 
\begin{equation}
 q_t(y,y') = \det\{a_{i,j}(y_i,y_j')\}_{1\leq i,j \leq n} \label{deta}
\end{equation}
where $a_{i,j}$ is given by 
\begin{equation}\label{aij}
 a_{i,j}(y,y') = (-1)^{i-1} \phi_t^{(j-i)}(y+y') 
                +(-1)^{i+j} \phi_t^{(j-i)}(y-y').
\end{equation}
\end{prop}
\noindent
The same type of formula was first obtained for the totally asymmetric 
simple exclusion process by Sch\"utz \cite{Schuetz1997b}. The formula 
for the $Z$ process was given as a Proposition 8 in
\cite{Warren2007}, see also \cite{Sa1998}. 

\begin{proofOF}{Proposition \ref{GY}}
For a fixed $y^\prime$, define $G_t(y,t)$ to be (\ref{deta}) as a function of $y$ and $t$.  
We check that $G$ satisfies (i) the heat equation, 
(ii) the boundary conditions $\frac{\partial G}{\partial y_1}|_{y_1=0}=0$,
$\frac{\partial G}{\partial y_i}|_{y_i=y_{i-1}}=0, ~i=2,3,\ldots,n$
and (iii) the initial conditions $G(y,t=0)=\prod_{i=1}^n \delta(y_i-y_i')$.

(i) holds since $\phi_t^{(k)}(y)$ for each $k$ satisfies the 
heat equation. (ii) follows from the relations,
$\frac{\partial}{\partial y} a_{1j}(y,y')|_{y=0} = 
\phi_t^{(j)}(y')+(-1)^{j+1}\phi_t^{(j)}(-y')=0$ and 
$\frac{\partial}{\partial y} a_{ij}(y,y') = -a_{i-1,j}(y,y')$. 
For (iii) we notice that the first term in (\ref{aij}) goes to zero as 
$t\to 0$ for $y,y'>0$ and the statement for the remaining part 
is shown in Lemma 7 in \cite{Warren2007}.  
\end{proofOF}

For $n=2m$, resp. $n=2m-1$ we take $(X(t),t\geq 0)$ to be Dyson
Brownian motion of type $C$, resp. of type $D$. The transition 
semigroup $\bigl(P_t; t \geq 0\bigr)$ of this process is given by (\ref{Pt}). 

\begin{figure}
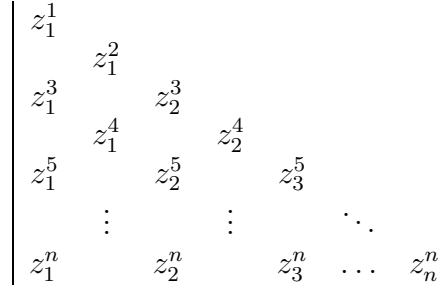

\begin{center}
\begin{tabular}{|*{9}{c}}
 $z_1^1$ &      &       &       &&      &       && \\
         &$z^2_1$&      &       &&      && \\
 $z^3_1$&       &$z^3_2$&&      &&& \\
         &$z^4_1$&      &$z^4_2$&       &&& \\
 $z^5_1$&       &$z^5_2$&       &$z^5_3$&& \\

 &$\vdots$&&    $\vdots$& &$\ddots$& \\
 $z^n_1$& & $z^n_2$& & $z^n_3$ &$\ldots$& $z^n_n$ \\
\end{tabular}
\end{center}
\caption{The set $\mathbb{K}$. The triangle represents 
the intertwining relations of the variables $z$ and the 
vertical line on the left indicates $z_1^{2k+1} \geq 0$,
see (\ref{inter1}),(\ref{inter2}).The set of variables on the bottom 
line is denoted by $b(z)$ and the one on the upper right line by $e(z)$.}
\label{sGZ}
\end{figure}

Let $\mathbb{K}$ denote the set with $n$ layers $z=(z^1,z^2,\ldots,z^n)$
where $z^{2k}=(z_1^{2k},z_2^{2k},\ldots,z_k^{2k})\in \mathbb{R}_+^k$,
$z^{2k-1}=(z_1^{2k-1},z_2^{2k-1},\ldots,z_k^{2k-1})\in \mathbb{R}_+^k$
and the intertwining relations,
\begin{equation}
 z_1^{2k-1} \leq z_1^{2k} \leq z_2^{2k-1} \leq z_2^{2k} \leq \ldots 
 \leq z_k^{2k-1} \leq z_k^{2k}
\label{inter1}
\end{equation}
and 
\begin{equation}
 0 \leq z_1^{2k+1} \leq z_1^{2k} \leq z_2^{2k+1} \leq z_2^{2k} \leq
 \ldots \leq z_k^{2k} \leq z_{k+1}^{2k+1}
\label{inter2}
\end{equation}
hold (Fig. \ref{sGZ}). 
Let $n=2m$ or $n=2m-1$ for some integer $m$. We define a kernel
$L^0$ from $E=\{0\leq x_1 \leq \ldots \leq x_m\}$ 
to $F=\{0\leq y_1 \leq \ldots \leq y_n\}$. For 
$z\in \mathbb{K}$, define $b(z) = z^n = (z_1^n,\ldots,z_m^n)\in E,$ 
$e(z) = (z_1^1,z_1^2,z_2^3,z_2^4,\ldots,z_m^n) \in F$ and 
$\mathbb{K}(x) = \{z\in\mathbb{K};b(z)=x\in E\},
\mathbb{K}[y] = \{z\in\mathbb{K};e(z)=y\in F\}$. 
The kernel $L^0$ is defined by 
\begin{equation}
 L^0 g (x)= \int_F L^0(x,dy)g(y)
 =
 \int_{\mathbb{K}(x)} g(e(z))dz.
\end{equation}
where the integrals are taken with respect to Lebesgue measure
but integrations with respect to $z$ on the RHS is for $b(z)=x$ fixed.

The function $h$ defined at \eqref{h} is equal to the Euclidean
volume of ${\mathbb K}(x)$.  Consequently we may define  $L$ to  be the Markov kernel $L(x,dy)=L^0(x,dy)/h(x)$.
In the remaining part of this section we show
\begin{prop}\label{LQPL}
\begin{equation}\label{LQPLe}
 L Q_t = P_t L.
\end{equation}
\end{prop}
\noindent
Now if we apply Lemma \ref{RP} with $f(x)=x_m$, $g(y)=y_n$ and 
the initial conditions starting from the origin     
we obtain (\ref{YX}).  

\begin{proofOF}{Proposition \ref{LQPL}} 
The kernels $P_t(x, \cdot)$ and $L(x, \cdot)$ are continuous in $x$.
Thus we may consider  $x$ in the interior of $E$, and  it is  enough to prove
\begin{equation} \label{LQPL0}
 (L^0 Q_t)(x,dy) = (P_t^0 L^0)(x,dy).
\end{equation}
From the definition of the kernel $L^0$, this is equivalent to 
showing 
\begin{equation}
 \int_{\mathbb{K}(x)} q_t(e(z),y)dz
 =
 \int_{\mathbb{K}[y]}p_t^0(x,b(z))dz
\label{KQKPint}
\end{equation}
where $q_t$ and $p^0$ are densities corresponding to $Q_t$ and
$P^0_t$. Integrations with respect to $z$ are on the LHS with 
$b(z)=x$ fixed and on the RHS with $e(z)=y$ fixed. 

Let us  consider the case where $n=2m$. 
Using the determinantal expressions for $q_t$ and $p_t^0$
we show that both sides of (\ref{KQKPint}) are equal to 
the determinant of size $2m$ whose $(i,j)$ matrix element is 
$a_{2i,j}(0,y_j)$ for $1\leq i \leq m, 1\leq j\leq 2m$ and
$a_{2m,j}(x_{i-m},y_j)$ for $m+1 \leq i \leq 2m, 1\leq j \leq 2m$.

The integrand of the LHS of (\ref{KQKPint}) is 
\begin{equation}
 q_t(e(z),y) = \det\{a_{i,j}(e(z)_i,y_j)\}_{1\leq i,j \leq 2m}
\end{equation}
with $b(z)=x$. We perform the integral with respect to 
$z^1,\ldots, z^{2m-1}$ in this order. 
After the integral up to $z^{2l-1}, 1\leq l \leq m$, we get the 
determinant of size $2m$ whose $(i,j)$ matrix element is 
$a_{2i,j}(0,y_j)$ for $1\leq i\leq l$, $a_{2l,j}(z_{i-l}^{2l},y_j)$
for $l+1\leq i\leq 2l$ and $a_{i,j}(e(z)_i,y_j)$ for 
$2l+1 \leq i \leq 2m$. Here we use a property of $a_{i,j}$, 
\begin{equation}
 a_{i,j}(y,y') = \int_y^{\infty} a_{i-1,j}(u,y')du, 
\end{equation}
and do some row operations in the determinant. The case for 
$l=m$ gives the desired expression.  

The integrand of the RHS of (\ref{KQKPint}) is
\begin{equation}
 p_t^0(x,z^{2m}) = \det(a_{2m,2m}(x_i,z_j^{2m}))_{1\leq i,j\leq m}
\end{equation}
with the condition $e(z)=y$. 
We perform the integrals with respect to 
$(z_1^{2m},\ldots,z_{m-1}^{2m}),(z_1^{2m-1},\ldots,z_{m-1}^{2m-1}),
\ldots, z_1^4,z_1^3$ in this order. 
We use properties of $a_{i,j}$,
\begin{align}
 a_{i,j}(y,y')  &= -\int_{y'}^{\infty} a_{i,j+1}(y,u)du , \\
a_{2i,2j}(x,0)=0, \;\;\; &a_{2i,2i-1}(0,y)=1,\;\;\; a_{2i , j }(0,y)=0,~ 2i \leq j. 
\end{align}
After each
integration corresponding to a layer of ${\mathbb K}$ we
simplify the determinant using column operations. We also   
expand the size of the determinant after an
integration corresponding to 
$(z_1^{2l},\ldots,z_{l-1}^{2l})$ for $1\leq l\leq m$, 
 by adding a new first row
\begin{multline}
\bigl(\underbrace{ 1,1,\ldots,1}_{l} ,\underbrace{ 0,0, \ldots
  ,0}_{2m-2l+1}\bigr) = \\
 \bigl( a_{2l,2l-1}(0,z^{2l-1}_1), \ldots,
a_{2l,2l-1}(0,z^{2l-1}_l),a_{2l,2l}(0,e(z)_{2l}), \ldots, a_{2l,2m}(0,e(z)_{2m}) )
\bigr)
\end{multline}
together with a new column.
After the integrals up to
$(z_1^{2l-1},\ldots,z_{l-1}^{2l-1})$ have been performed,  we obtain the determinant of size $2m-l+1$,
\begin{equation}
\begin{vmatrix}
  a_{2(l+i-1),2(l-1)}(0,z_j^{2(l-1)}) & a_{2(l+i-1),j+l-1}(0,e(z)_{j+l-1}) \\
  a_{2m,2(l-1)}(x_{i-m+l-1},z_j^{2(l-1)}) 
  & a_{2m,j+l-1}(x_{i-m+l-1},e(z)_{j+l-1}) 
\end{vmatrix}. 
\end{equation}
Here $1\leq i \leq m-l+1$ (resp. $m-l+2 \leq i \leq 2m-l+1$) 
for the upper expression (resp. the lower expression) and 
$1\leq j\leq l-1$ (resp. $l\leq j \leq 2m-l+1$) for the 
left (resp.  right) expression. 
For $l=1$ this reduces to the same determinant as for the 
LHS.  

The case $n=2m-1$ is almost identical. Similar arguments show that both
sides of \eqref{LQPL0} are equal to a determinant size $2m-1$ whose $(i,j)$ matrix element is 
$a_{2i,j}(0,y_j)$ for $1\leq i \leq m-1, 1\leq j\leq 2m-1$ and
$a_{2m-1,j}(x_{i-m+1},y_j)$ for $m+1 \leq i \leq 2m-1, 1\leq j \leq 2m-1$.
\end{proofOF}

\section{Proof of proposition \ref{maxWWr}}{\label{WWr}}
Using (\ref{YinB}) repeatedly, one has
\begin{equation}
 Y_n(t) 
 = 
 \sup_{0\leq t_1 \leq \ldots \leq t_n\leq t}
 \sum_{i=1}^n (B_i(t_{i+1})-B_i(t_i))
\end{equation}
with $t_{n+1}=t$. By renaming $t-t_{n-i+1}$ by $t_i$ and changing 
the order of the summation, we have
\begin{equation}
 Y_n(t) 
 = 
 \sup_{0\leq t_1 \leq \ldots \leq t_n\leq t}
 \sum_{i=1}^n (B_{n-i+1}(t-t_{i+1})-B_{n-i+1}(t-t_i)).
\end{equation}
Since $\tilde{B}_i(s):=B_{n-i+1}(t)-B_{n-i+1}(t-s) \overset{d}{=}
B_i(s)$, 
\begin{equation}
 Y_n(t) 
 \overset{d}{=} 
 \sup_{0\leq t_1 \leq \ldots \leq t_n\leq t}
 \sum_{i=1}^n (B_i(t_i)-B_i(t-t_{i-1}))
 =
 \sup_{0\leq s\leq t} Z_n(t).
\end{equation}

\providecommand{\bysame}{\leavevmode\hbox to3em{\hrulefill}\thinspace}
\providecommand{\MR}{\relax\ifhmode\unskip\space\fi MR }
\providecommand{\MRhref}[2]{%
  \href{http://www.ams.org/mathscinet-getitem?mr=#1}{#2}
}
\providecommand{\href}[2]{#2}

\end{document}